\documentclass[11pt,a4paper]{amsart}
\usepackage{amssymb}
\usepackage[a4paper,centering,twoside,textwidth=6in]{geometry}

\newcommand{\De}{\Delta}
\newcommand{\al}{\alpha}
\newcommand{\be}{\beta}
\newcommand{\ga}{\gamma}

\newcommand{\eps}{\epsilon}

\DeclareMathOperator{\PSL}{PSL} \DeclareMathOperator{\PGL}{PGL}
\DeclareMathOperator{\PSU}{PSU} \DeclareMathOperator{\SL}{SL}
 
 \DeclareMathOperator{\Aut}{Aut}
\DeclareMathOperator{\tr}{tr} 
\DeclareMathOperator{\Irr}{Irr}

\newcommand{\bZ}{\mathbb{Z}}
\newcommand{\BF}{\mathbb{F}}

\theoremstyle{plain}
\newtheorem{prop}{Proposition}%[section]
\newtheorem{lemma}[prop]{Lemma}

\newtheorem{conj}[prop]{Conjecture}
\newtheorem{thm}[prop]{Theorem}

\theoremstyle{definition}

\newtheorem{defn}[prop]{Definition}
\newtheorem{quest}[prop]{Question}

\title{Beauville surfaces and probabilistic group theory}
\author{Shelly Garion}
\thanks{The author was supported by the SFB 878 ``Groups, Geometry and Actions''.}

\address{Fachbereich Mathematik und Informatik, Universit\"at
M\"unster, Einsteinstrasse 62, D-48149 M\"unster, Germany}
\email{shelly.garion@uni-muenster.de}

\begin{document}

\maketitle

%%%%%%%%%%%%%%%%%%%%%%%%%%%%%%%%%%%%%%%%%%%%%%%%%%%%%%%
\begin{abstract}
A Beauville surface is a complex algebraic surface that can be
presented as a quotient of a product of two curves by a suitable
action of a finite group. Bauer, Catanese and Grunewald have been
able to intrinsically characterize the groups appearing in minimal
presentations of Beauville surfaces in terms of the existence of a
so-called ''Beauville structure``. They conjectured that all finite
simple groups, except $A_5$, admit such a structure. This conjecture
has recently been proved by Guralnick-Malle and
Fairbairn-Magaard-Parker.

In this survey we demonstrate another approach towards the proof of
this conjecture, based on probabilistic group-theoretical methods,
by describing the following three works. The first is the work of
Garion, Larsen and Lubotzky, showing that the above conjecture holds
for almost all finite simple groups of Lie type. The second is the
work of Garion and Penegini on Beauville structures of alternating
groups, based on results of Liebeck and Shalev, and the third is the
case of the group $\PSL_2(p^e)$, in which we give bounds on the probability
of generating a Beauville structure. We also discuss other related
problems regarding finite simple quotients of hyperbolic triangle
groups and present some open questions and conjectures.
\end{abstract}

%%%%%%%%%%%%%%%%%%%%%%%%%%%%%%%%%%%%%%%%%%%%%%%%%%%%%
\section{Beauville surfaces and Beauville structures}
%%%%%%%%%%%%%%%%%%%%%%%%%%%%%%%%%%%%%%%%%%%%%%%%%%%%%

A \emph{Beauville surface} $S$ (over $\mathbb{C}$) is a particular
kind of surface isogenous to a \emph{higher product of curves},
i.e., $S=(C_1 \times C_2)/G$ is a quotient of a product of two
smooth curves $C_1$ and $C_2$ of genus at least two, modulo a free
action of a finite group $G$ which acts faithfully on each curve.
For Beauville surfaces the quotients $C_i/G$ are isomorphic to
$\mathbb{P}^1$ and both projections $C_i \rightarrow C_i/G \cong
\mathbb{P}^1$ are coverings branched over three points. A Beauville
surface is in particular a minimal surface of general type.

Beauville~\cite{Be} constructed a minimal surface of general type
$S$ with $K^2_S=8$ and $p_g=q=0$ in the following way: take two
curves $C_1=C_2$ given by the Fermat equation $x^5+y^5+z^5=0$ and
$G$ the group $(\bZ/5\bZ)^2$ acting on $C_1 \times C_2$ by
\[
 (a,b) \cdot ([x:y:z],[u:v:w]) =
 ([\xi^a x: \xi^b y : z],[\xi^{a+3b}u : \xi^{2a+4b}v : w]),
\]
where $\xi = e^{\frac{2\pi i}{5}}$ and $a,b \in \bZ/5\bZ$. Then
define $S$ by the quotient $(C_1 \times C_2)/G$. Moreover $C_i
\rightarrow C_i/G \cong \mathbb{P}^1_{\mathbb{C}}$ and both covers
are branched in exactly three points. Curves with such properties
are said to be \emph{triangle curves}.

Inspired by this construction Catanese~\cite{Cat00} observed that in
general if $C_1$ and $C_2$ are two triangle curves with group $G$,
if the action of $G$ on the product $C_1 \times C_2$ is free, then
$S=(C_1 \times C_2)/G$ is a \emph{strongly rigid surface}, i.e., if
$S'$ is another surface homotopically equivalent to $S$ then $S'$ is
either biholomorphic  or antibiholomorphic to $S$. He proposed to
name these surfaces \emph{Beauville surfaces}.

A Beauville surface $S$ is either of \emph{mixed} or \emph{unmixed}
type according respectively as the action of $G$ exchanges the two
factors (and then $C_1$ and $C_2$ are isomorphic) or $G$ acts
diagonally on the product $C_1 \times C_2$. The subgroup $G_0$ (of
index $\leq 2$) of $G$ which preserves the ordered pair $(C_1,C_2)$
is then respectively of index $2$ or $1$ in $G$. Any Beauville
surface $S$ can be presented in such a way that the subgroup $G_0$
of $G$ acts effectively on each of the factors $C_1$ and $C_2$.
Catanese called such a presentation \emph{minimal} and proved its
uniqueness in \cite{Cat00}. In this survey we shall consider only
Beauville surfaces of unmixed type so that $G_0 = G$.

An extensive research on Beauville surfaces was initiated by the
collaboration of Bauer, Catanese and Grunewald \cite{BCG05,BCG06}.
They have been able to intrinsically characterize the groups
appearing in minimal presentations of unmixed Beauville surfaces in
terms of the existence of the so-called unmixed \emph{Beauville
structure}.

\begin{defn}\label{defn.beau}
An \emph{unmixed Beauville structure} for a finite group $G$ is a
quadruple $(x_1, y_1; x_2, y_2)$ of elements of $G$, which
determines two triples $(x_1,y_1,z_1)$ and $(x_2,y_2,z_2)$
satisfying:
\begin{enumerate}\renewcommand{\theenumi}{\it \roman{enumi}}
    \item $x_1y_1z_1=1$ and $x_2y_2z_2=1$,
    \item $\langle x_1,y_1 \rangle = G$ and $\langle x_2,y_2 \rangle=G$,
    \item $\Sigma(x_1,y_1,z_1) \cap \Sigma(x_2,y_2,z_2)= \{1\}$, \\
    where $\Sigma(x,y,z)$ is the union of the conjugacy classes of
    all powers of $x$, all powers of $y$, and all powers of
    $z$.
\end{enumerate}
\end{defn}

Moreover denoting the order of an element $g$ in $G$ by $|g|$, we
define the \emph{type} $\tau$ of $(x,y,z)$ to be the triple
$(|x|,|y|,|z|)$. In this situation, we say that $G$ admits an
\emph{unmixed Beauville structure of type} $(\tau_1,\tau_2)$.

The question whether a finite group admits an unmixed Beauville
structure of a given type is closely related to the question whether
it is a quotient of certain triangle groups. More precisely, a
necessary condition for a finite group $G$ to admit an unmixed
Beauville structure of type $(\tau_1,\tau_2) =
\bigl((r_1,s_1,t_1),(r_2,s_2,t_2)\bigr)$ is that $G$ is a quotient
with torsion free-kernel of the triangle groups $T_{r_1,s_1,t_1}$
and $T_{r_2,s_2,t_2}$, where
\[
    T_{r,s,t} = \langle x,y,z: x^{r}=y^{s}=z^{t}=xyz=1 \rangle.
\]
Indeed, conditions \emph{(i)} and \emph{(ii)} of
Definition~\ref{defn.beau} are equivalent to the condition that $G$
is a quotient of each of the triangle groups
$T_{|x_i|,|y_i|,|z_i|}$, for $i\in \{1,2\}$, with torsion-free
kernel.

When investigating the existence of an unmixed Beauville structure
for a finite group, one can consider only types $(\tau_1,\tau_2)$,
where for $i\in \{1,2\}$, $\tau_i = (r_i,s_i,t_i)$ satisfies $1/r_i
+ 1/s_i + 1/t_i < 1$. Then $T_{r_i,s_i,t_i}$ is a (infinite
non-soluble) hyperbolic triangle group and we say that $\tau_i$ is
\emph{hyperbolic}. Indeed, if $1/r_i + 1/s_i + 1/t_i > 1$ then
$T_{r_i,s_i,t_i}$ is a finite group, and moreover, it is either
dihedral or isomorphic to one of $A_4$, $A_5$ or $S_4$. By
\cite[Proposition 3.6 and Lemma 3.7]{BCG05}, in these cases $G$
cannot admit an unmixed Beauville structure. If $1/r_i + 1/s_i +
1/t_i = 1$ then $T_{r_i,s_i,t_i}$ is one of the (soluble infinite)
``wall-paper'' groups, and by \cite[\S 6]{BCG05}, none of its finite
quotients can admit an unmixed Beauville structure.

Observe that condition~\textit{(iii)} of Definition~\ref{defn.beau}
is clearly satisfied under the assumption that $r_1 s_1 t_1$ is
coprime to $r_2 s_2 t_2$. However this assumption is not always
necessary, as demonstrated by many examples, such as Beauville's
original construction, abelian groups~\cite[Theorem 3.4]{BCG05},
alternating groups~\cite[Theorem 1.2]{GP} and the group
$\PSL_2(p^e)$~\cite{Ga}.

%%%%%%%%%%%%%%%%%%%%%%%%%%%%%%%%%%%%%%%%%%%%%%%%%%%%%
\section{Beauville surfaces and finite simple groups}
%%%%%%%%%%%%%%%%%%%%%%%%%%%%%%%%%%%%%%%%%%%%%%%%%%%%%

A considerable effort has been made to classify the finite simple
groups which admit an unmixed Beauville structure. We recall that by
\emph{the classification theorem of finite simple groups}, any
finite simple group belongs to one of the following families: the
cyclic groups $Z_p$ of prime order; the alternating permutation
groups $A_n (n \geq 5)$; the finite simple groups of Lie type,
defined over finite fields (e.g. $\PSL_n(q))$; and finally the 26
so-called sporadic groups.

A finite abelian simple group clearly does not admit an unmixed
Beauville structure as given a prime $p$, any pair $(a,b)$ of
elements of the cyclic group $Z_p$ of prime order $p$ generating it
satisfies $\Sigma(a,b,c) = Z_p$. In fact Bauer, Catanese and
Grunewald showed in \cite[Theorem 3.4]{BCG05} that the only finite
abelian groups admitting an unmixed Beauville structure are the
abelian groups of the form $Z_n \times Z_n$ where $n$ is a positive
integer coprime to $6$. (Here $Z_n$ denotes a cyclic group of order
$n$.)

In \cite{BCG05} the authors also provide the first results on finite
non-abelian simple groups admitting an unmixed Beauville structure.
More precisely they show that the alternating groups of sufficiently
large order admit an unmixed Beauville structure, as well as the
projective special linear groups $\PSL_2(p)$ where $p>5$ is a prime.
Moreover using computational methods, they checked that every finite
non-abelian simple group of order less than $50000$ admits an
unmixed Beauville structure with the exception of the alternating
group $A_5$. Based on these results and the latter observation, they
conjectured that all finite non-abelian simple groups admit an
unmixed Beauville structure with the exception of $A_5$.

This conjecture has received much attention and has recently been
proved to hold. Concerning the simple alternating groups, it was
established in \cite{FG} that $A_5$ is indeed the only one not
admitting an unmixed Beauville structure. In \cite{FJ,GP}, the
conjecture is shown to hold for the projective special linear groups
$\PSL_2(q)$ (where $q > 5$), the Suzuki groups $^2B_2(q)$ and the
Ree groups $^2G_2(q)$ as well as other families of finite simple
groups of Lie type of small rank. More precisely, the projective
special and unitary groups $\PSL_3(q)$, $\PSU_3(q)$, the simple
groups $G_2(q)$ and the Steinberg triality groups $^3D_4(q)$ are
shown to admit an unmixed Beauville structure if $q$ is large (and
the characteristic $p$ is greater than $3$ for the simple
exceptional groups of type $G_2$ or $^3D_4$). The next major result
concerning the investigation of the conjecture with respect to the
finite simple groups of Lie type was pursued by Garion, Larsen and
Lubotzky who showed in \cite{GLL} that the conjecture holds for
finite non-abelian simple groups of sufficiently large order. The
final step regarding the investigation of the conjecture was carried
out by Guralnick and Malle \cite{GM} and Fairbairn, Magaard and
Parker \cite{FMP} who established its veracity in general, namely,

\begin{thm}\cite{GM}.\label{thm.Beau.FSG}
Any finite non-abelian simple group, except $A_5$, admits an unmixed
Beauville structure.
\end{thm}

There has also been an effort to classify the finite quasisimple
groups and almost simple groups which admit an unmixed Beauville
structure. Recall that a finite group $G$ is \emph{quasisimple}
provided $G/Z(G)$ is a non-abelian simple group and $G = [G,G]$. In
\cite{FJ} it was shown that $\SL_2(q)$ (for $q>5$) admits an unmixed
Beauville structure. Fairbairn, Magaard and Parker \cite{FMP} proved
the following general result.

\begin{thm}\cite{FMP}.
With the exceptions of $\SL_2(5)$ and $\PSL_2(5) \cong \SL_2(4)
\cong A_5$, every finite quasisimple group admits an unmixed
Beauville structure.
\end{thm}

Recall that a group $G$ is called \emph{almost simple} if there is a
non-abelian simple group $G_0$ such that $G_0 \leq G \leq
\Aut(G_0)$. By \cite{BCG06,FG} the symmetric groups $S_n$ (where
$n\geq 5$) admit an unmixed Beauville structure, and by~\cite{Ga}
the group $\PGL_2(p^e)$ admits such a structure.

Moreover for the alternating and symmetric groups Garion and
Penegini~\cite{GP} proved another conjecture that Bauer, Catanese
and Grunewald proposed in~\cite{BCG05}, that almost all of these
groups admit a Beauville structure with fixed \emph{type}, namely,

\begin{thm}\cite[Theorem 1.2]{GP}.\label{thm.Beau.An}
If $\tau_1=(r_1,s_1,t_1)$ and $\tau_2=(r_2,s_2,t_2)$ are two
hyperbolic types, then almost all alternating groups $A_n$ admit an
unmixed Beauville structure of type $(\tau_1,\tau_2)$.
\end{thm}

A similar theorem also applies for symmetric groups, see~\cite{GP},
and a similar conjecture was raised in~\cite{GP}, replacing $A_n$ by
a finite simple classical group of Lie type of sufficiently large
Lie rank, namely,

\begin{conj}\cite[Conjecture 1.7]{GP}.
Let $\tau_1=(r_1,s_1,t_1)$ and $\tau_2=(r_2,s_2,t_2)$ be two
hyperbolic types. If $G$ is a finite simple classical group of Lie
type of Lie rank large enough, then it admits an unmixed Beauville
structure of type $(\tau_1,\tau_2)$.
\end{conj}

In contrast, when the Lie rank is very small, as in the case of
$\PSL_2(q)$, such a conjecture does not hold, as demonstrated
in~\cite{Ga}, where there is a characterization of the types of
Beauville structures for these groups.

\medskip

It is well known that almost all pairs of elements in a finite
simple (non-abelian) group are generating pairs~\cite{Dix,KL,LS95},
hence the following question was raised in the workshop ``Beauville
surfaces and groups 2012''.

\begin{quest}
Let $G$ be a finite (non-abelian) simple group. What is the
probability $P(G)$ that for four random elements $x_1,y_1,x_2,y_2
\in G$ the quadruple $(x_1,y_1;x_2,z_2)$ is an unmixed Beauville
structure for $G$?

In particular, is it true that if $G=A_n$ or $G=G_n(q)$, a finite
simple group of Lie type of Lie rank $n$, then $P(G) \rightarrow 1$
as $n \rightarrow \infty$?
\end{quest}

Two interesting comments were made during the workshop regarding
this question. The first comment, due to Malle, is that for finite
simple groups of Lie type of bounded Lie rank, $P(G)$ does not go to
$1$, and it is bounded above by a function of the rank. The second
comment, due to Magaard, is that the techniques in~\cite{FMP}
demonstrate that one can generate many unmixed Beauville structures
for the finite simple groups of Lie type, allowing to obtain a
constant lower bound on $P(G)$, when $G$ is a finite simple
classical group. In the specific case where $G=\PSL_2(q)$ we give
the following bounds on the probability of $P(G)$
(see~\S\ref{sect.prob.PSL} for the proof).

\begin{thm}\label{thm.prob.PSL}
Let $G=\PSL_2(q)$.
\begin{itemize}
\item If $q$ is odd then $\frac{1}{32} - \eps_q \leq P(G) \leq
\frac{15}{16} + \eps_q$,
\item If $q$ is even then $\frac{1}{32} -
\eps_q \leq P(G) \leq \frac{35}{36} + \eps_q$,
\end{itemize}
where $\eps_q \rightarrow 0$ as $q \rightarrow \infty$.
\end{thm}

%%%%%%%%%%%%%%%%%%%%%%%%%%%%%%%%%%%%%%%%%%%%%%%%%%%%%%%%%%%%%%%%%%
\section{Hyperbolic triangle groups and their finite quotients}
%%%%%%%%%%%%%%%%%%%%%%%%%%%%%%%%%%%%%%%%%%%%%%%%%%%%%%%%%%%%%%%%%%
Since for a finite group $G$ which admits an unmixed Beauville
structure there exists an epimorphism from a hyperbolic triangle
group to $G$, we recall in this section some results on finite
quotients of hyperbolic triangle groups.

A \emph{hyperbolic triangle group} $T$ is a group with presentation
\[
 T=T_{r,s,t} = \langle x,y,z: x^r=y^s=z^t=xyz=1 \rangle,
\]
where $(r,s,t)$ is a triple of positive integers satisfying the
condition $1/r+1/s+1/t < 1$. Geometrically, let $\De$ be a
hyperbolic triangle group having angles of sizes $\pi/r$, $\pi/s$,
$\pi/t$, then $T$ can be viewed as the group generated by rotations
of angles $\pi/r$, $\pi/s$, $\pi/t$ around the corresponding
vertices of $\De$ in the hyperbolic plane $\mathbb{H}^2$. Moreover,
a hyperbolic triangle group $T_{r,s,t}$ has positive measure
$\mu(T_{r,s,t})$ where $\mu(T_{r,s,t})=1-(1/r+1/s+1/t)$. As
hyperbolic triangle groups are infinite and non-soluble it is
interesting to study their finite quotients, particularly the simple
ones.

A hyperbolic triangle group $T_{r,s,t}$ has minimal measure when
$(r,s,t)=(2,3,7)$. The group $T_{2,3,7}$ is also called the
$(2,3,7)$-triangle group and its finite quotients are also known as
\emph{Hurwitz groups}. These are named after Hurwitz who showed in
the late nineteenth century that if $S$ is a compact Riemann surface
of genus $h\geq 2$ then $|\Aut S| \leq 84(h-1)$ and this bound is
attained if and only if $\Aut S$ is a quotient of the triangle group
$T_{2,3,7}$. Following this result, much effort has been given to
classify Hurwitz groups, especially the simple ones, see for example
\cite{Co90} for a historical survey, and \cite{Co10,TVsum} for the
current state of the art.

Most alternating groups are Hurwitz as shown by Conder (following
Higman) who proved in \cite{Co80} that if $n > 167$ then the
alternating group $A_n$ is a quotient of $T_{2,3,7}$. Concerning the
finite simple groups of Lie type, there is a dichotomy with respect
to their occurrence as quotients of $T_{2,3,7}$ depending on whether
the Lie rank is large or not. Indeed as shown in \cite{LT} many
classical groups of large rank are Hurwitz (and there is no known
example of classical groups of large rank which are not Hurwitz). As
an illustration by \cite{LTW} if $n \geq 267$ then the projective
special linear group $\PSL_n(q)$ is Hurwitz for any prime power $q$.
The behavior of finite simple groups of Lie type of relatively low
rank with respect to the Hurwitz generation problem is rather
sporadic. As an illustration by respective results of
\cite{Coh,TVsum,Mal,Mac}, $\PSL_3(q)$ is Hurwitz if and only if
$q=2$, $\PSL_4(q)$ is never Hurwitz, $G_2(q)$ is Hurwitz for $q\geq
5$, and $\PSL_2(p^e)$ is Hurwitz if and only if $e = 1$ and $p
\equiv 0,\pm 1 \bmod 7$, or $e=3$ and $p \equiv \pm 2,\pm 3 \bmod
7$. Therefore, unlike the alternating groups, there are finite
simple groups of Lie type of large order which are not quotients of
$T_{2,3,7}$. As for the 26 sporadic finite simple groups, 12 of them
are Hurwitz (including the Monster~\cite{Wi}) while the other 14
groups are not.

\medskip

Turning to general hyperbolic triples $(r,s,t)$ of integers, Higman
had already conjectured in the late 1960s that every hyperbolic
triangle group has all but finitely many alternating groups as
quotients. This was eventually proved by Everitt~\cite{Ev}, namely,

\begin{thm}\cite{Ev}.
For any hyperbolic triangle group $T=T_{r,s,t}$, if $n \geq
n_0(r,s,t)$ then the alternating group $A_n$ is a quotient of $T$.
\end{thm}

Later, Liebeck and Shalev~\cite{LS04} gave an alternative proof to
Higman's Conjecture based on probabilistic group theory, and
moreover they have conjectured the following.

\begin{conj}\cite{LS05F}.
For any hyperbolic triangle group $T=T_{r,s,t}$, if $G=G_n(q)$ is a
finite simple classical group of Lie rank $n \geq n_0(r,s,t)$, then
the probability that a randomly chosen homomorphism from $T$ to $G$
is an epimorphism tends to $1$ as $|G| \rightarrow \infty$.
\end{conj}
% new results of Bezrukavnikov, Liebeck, Shalev ??? %

This conjecture has been proved by Marion~\cite{Mar11,MarPre} 
for certain families of groups of small Lie rank and certain 
triples $(r,s,t)$. For example, take $(r,s,t)$ to be a
hyperbolic triple of odd primes and $G = \PSL_3(q)$ or $\PSU_3(q)$
containing elements of orders $r, s$ and $t$, then the conjecture
holds. As another example, if $(r,s,t)$ is a hyperbolic triple of
primes and $G =\ ^2B_2(q)$ or $^2G_2(q)$ contains elements of orders
$r, s$ and $t$, then the conjecture also holds.

However, for finite simple groups of small Lie rank such a
conjecture does not hold \emph{in general}, and it fails to hold in
the case of $\PSL_2(q)$. Indeed, Langer and Rosenberger~\cite{LR1}
and Levin and Rosenberger~\cite{LR2} had generalized the
aforementioned result of Macbeath, and determined, for a given prime
power $q=p^e$, all the triples $(r,s,t)$ such that $\PSL_2(q)$ is a
quotient of $T_{r,s,t}$, with torsion-free kernel. It follows that
if $(r,s,t)$ is hyperbolic, then for almost all primes $p$, there is
precisely one group of the form $\PSL_2(p^e)$ or $\PGL_2(p^e)$ which
is a homomorphic image of $T_{r,s,t}$ with torsion-free kernel. We
note that this result can also be obtained by using other
techniques. Firstly, Marion~\cite{Mar09} has recently provided a
proof for the case where $r,s,t$ are primes relying on probabilistic
group theoretical methods. Secondly, it also follows from the
representation theoretic arguments of Vincent and
Zalesski~\cite[Theorems 2.9 and 2.11]{VZ}. Such methods can be used
for dealing with other families of finite simple groups of Lie type,
see for example~\cite{Mar10,MTZ,TV,TZ,VZ}.

Recently, a new approach was presented by Larsen, Lubotzky and
Marion~\cite{LLM}, based on the theory of representation varieties
(via deformation theory). They prove a conjecture of
Marion~\cite{Mar10} showing that various finite simple groups are
\emph{not} quotients of $T_{r,s,t}$, as well as positive results
showing that many finite simple groups are quotients of $T_{r,s,t}$.

%%%%%%%%%%%%%%%%%%%%%%%%%%%%%%%%%%%%%%%%%%%%%%%%%%%%%%%%%%%
\section{Beauville structures for the group $\PSL_2(q)$}\label{sect.PSL2}
%%%%%%%%%%%%%%%%%%%%%%%%%%%%%%%%%%%%%%%%%%%%%%%%%%%%%%%%%%%
In this section we discuss the specific case of $\PSL_2(q)$, and
briefly sketch the proof of Garion and Penegini~\cite{GP} for the
following theorem, which is based on results of Macbeath~\cite{Mac}.

\begin{thm}\label{thm.PSL.Beau}\cite{FJ,GP}.
Let $p$ be a prime number, and assume that $q = p^e$ is at least
$7$. Then the group $\PSL_2(q)$ admits an unmixed Beauville
structure.
\end{thm}

In addition, we bound the probability that four random elements
in $\PSL_2(q)$ generate an unmixed Beauville structure and prove
Theorem~\ref{thm.prob.PSL}.

%%%%%%%%%%%%%%%%%%%%%%%%%%%%%%%%%%%%%%%%%%%%%%%%%%%%%%%%%%%%%%
\subsection{Sketch of the proof of Theorem~\ref{thm.PSL.Beau}}
%%%%%%%%%%%%%%%%%%%%%%%%%%%%%%%%%%%%%%%%%%%%%%%%%%%%%%%%%%%%%%
In order to construct an unmixed Beauville structure for $\PSL_2(q)$
one needs to find a quadruple $(A_1, B_1; A_2, B_2)$ of elements of
$\PSL_2(q)$ satisfying the three conditions given in
Definition~\ref{defn.beau}. This can be done directly by finding
specific elements in the group satisfying these conditions
(see~\cite{FJ}) or indirectly by using the following results of
Macbeath~\cite{Mac} (see~\cite{Ga,GP}).

\begin{thm}\cite[Theorem 1]{Mac}.
For every $\al,\be,\ga \in \BF_q$ there exist three matrices $A,B,C
\in \SL_2(q)$ satisfying $\tr(A)=\al$, $\tr(B)=\be$, $\tr(C)=\ga$
and $ABC=I$.
\end{thm}

This theorem immediately implies Condition $\textit{(i)}$.

Moreover, Macbeath~\cite{Mac} classified the pairs of elements in
$\PSL_2(q)$ in a way which makes it easy to decide what kind of
subgroup they generate. By~\cite[Theorem 2]{Mac}, a triple
$(\al,\be,\ga) \in \BF_q^3$ is \emph{singular}, namely
$\al^2+\be^2+\ga^2-\al\be\ga - 4 = 0$, if and only if for the
corresponding triple of matrices $(A, B,C)$, the group generated by
the images of $A$ and $B$ is a \emph{structural subgroup} of
$\PSL_2(q)$, that is a subgroup of the Borel or a cyclic subgroup.

%We briefly recall the subgroup structure of $\PSL_2(q)$
%(see~\cite{Di,Su}). The subgroups isomorphic to $\PSL_2(q_1)$ (or
%$\PGL_2(q_1)$), where $q=q_1^m \ (m \in \mathbb{N})$, are usually
%called \emph{subfield subgroups} (since $\BF_{q_1}$ is a subfield of
%$\BF_q$). Since $A_4$, $S_4$, $A_5$ and dihedral groups correspond
%to the finite triangle groups, that is, triangle groups $T_{r,s,t}$
%such that $1/r+1/s+1/t> 1$, we will call them \emph{small
%subgroups}. For convenience we will refer to the other subgroups,
%namely subgroups of the Borel and cyclic subgroups, as
%\emph{structural subgroups}.

Hence, in order to verify Condition $\textit{(ii)}$, one needs to
show that the subgroup generated by $A,B \in \PSL_2(q)$ is neither
a structural subgroup (using the aforementioned result of
Macbeath), not a dihedral subgroup, not one of the small subgroups 
$A_4$, $S_4$ or $A_5$, and not a subfield subgroup 
(namely, isomorphic to $\PSL_2(q_1)$ or to $\PGL_2(q_1)$, where $q=q_1^m$) hence it must be $\PSL_2(q)$ itself, as the subgroup structure of $\PSL_2(q)$
is well-known (see e.g.~\cite{Di,Su}).
For example, Condition $\textit{(ii)}$ is always satisfied
when $q\geq 13$ and the orders $|A|=|B|=|C|=(q-1)/d$ or $(q+1)/d$,
where $d=\gcd(2,q-1)$.

Condition~$\textit{(iii)}$ is clearly satisfied under the assumption
that the product of the orders $|A_1|\cdot|B_1|\cdot|C_1|$ is
coprime to $|A_2|\cdot|B_2|\cdot|C_2|$. For example, for any $q>7$
the group $\PSL_2(q)$ admits unmixed Beauville structures of types
\[
    \left(\bigl(\frac{q-1}{d},\frac{q-1}{d},\frac{q-1}{d}\bigr),
    \bigl(\frac{q+1}{d},\frac{q+1}{d},\frac{q+1}{d}\bigr)\right),
\]
and
\[
    \left(\bigl(\frac{q-1}{d},\frac{q-1}{d},\frac{q-1}{d}\bigr),
    \bigl(\frac{q+1}{d},\frac{q+1}{d},p\bigr)\right),
\]
appearing in~\cite{GP} and~\cite{FJ} respectively.

However this assumption is not always necessary. Indeed,
by~\cite{Ga}, $\PSL_2(q)$ (where $q=p^{2e}$, $p$ an odd prime)
always admits unmixed Beauville structures of types $((p, p, t_1),
(p, p, t_2))$ for certain $t_1, t_2$ dividing $(q-1)/2$, $(q+1)/2$
respectively.

This approach can be effectively used to construct many unmixed
Beauville structures for $\PSL_2(q)$, and in~\cite{Ga} there is a
characterization of the types of unmixed Beauville structures for
this group.

%%%%%%%%%%%%%%%%%%%%%%%%%%%%%%%%%%%%%%%%%%%%%%%%
\subsection{Proof of Theorem~\ref{thm.prob.PSL}}
\label{sect.prob.PSL}
%%%%%%%%%%%%%%%%%%%%%%%%%%%%%%%%%%%%%%%%%%%%%%%%
The proof relies on considering the various types of elements in
$G=\PSL_2(q)$. Recall that an element in $G$ is called \emph{split}
if its order divides $(q-1)/d$ (where $d=\gcd(2,q-1)$),
\emph{non-split} if its order divides $(q+1)/d$, and
\emph{unipotent} if its order is $p$ (and then the trace of its
pre-image in $\SL_2(q)$ equals $\pm 2$).

It is well-known that there are roughly $q^3$ matrices in
$\SL_2(q)$, and moreover, for any $\al \in \BF_q$ the number of
matrices $A \in \SL_2(q)$ with $\tr(A)=\al$ is roughly $q^2$ (see
for example \cite[Table 1]{BG}). In addition, the probability that a
random element in $\BF_q$ is a trace of a split (respectively,
non-split) matrix in $\SL_2(q)$ goes to $1/2$ as $q \rightarrow
\infty$ (see~\cite[Lemma 2]{Mac}). Therefore, probabilistically, the
number of unipotents in $G$ is negligible, and moreover, if we
denote by $P^s_q$ (respectively $P^n_q)$ the probability that a
random element in $G$ is split (respectively, non-split) then $P^s_q
\rightarrow \frac{1}{2}$ and $P^n_q \rightarrow \frac{1}{2}$ as $q
\rightarrow \infty$.

By~\cite[Proposition 7.2]{BG} it follows that for any non-singular
triple $(\al,\be,\ga)\in \BF_q^3$, the number of triples $(A,B,C)
\in \SL_2(q)^3$ satisfying $\tr(A)=\al$, $\tr(B)=\be$, $\tr(C)=\ga$
and $ABC=I$ is roughly $q^3$. By~\cite[Lemma 3]{Mac} almost all
triples in $\BF_q^3$ are non-singular. Since the probability that a
random triple of element in $\BF_q^3$ contains only traces of split
(respectively, non-split) matrices goes to $\frac{1}{8}$ as $q
\rightarrow \infty$, it follows that the probability that two random
elements $A,B \in G$ satisfy that $A,B$ and $AB$ are all split
(respectively, non-split), goes to $\frac{1}{8}$ as $q \rightarrow
\infty$.

In order to obtain a lower bound, observe that the
probability that two random elements $A,B \in G$ do {\bf not}
generate $G$ goes to $0$ as $q \rightarrow \infty$ 
(this can be deduced from the aformentioned results of Macbeath~\cite{Mac}, 
or alternatively, as a specific case of~\cite{KL}). Namely,
\[
    \frac{\#\{(A,B) \in G^2: \langle A,B \rangle \neq G\}}{|G|^2}
    \leq \eps'_q,
\]
where $\eps'_q \rightarrow 0$ as $q \rightarrow \infty$. Hence,
\[
    \frac{\#\{(A,B) \in G^2: A,B,AB \text{ are split and } \langle A,B \rangle = G\}}{|G|^2}
    \geq \frac{1}{8} - \eps'_q,
\]
and
\[
    \frac{\#\{(A,B) \in G^2: A,B,AB \text{ are non-split and } \langle A,B \rangle = G\}}{|G|^2}
    \geq \frac{1}{8} - \eps'_q.
\]
Since $(q-1)/d$ and $(q+1)/d$ are relatively prime then
Condition~$\textit{(iii)}$ is immediately satisfied when
$A_1,B_1,A_1B_1$ are all split and $A_2,B_2,A_2B_2$ are all
non-split (and vice-versa). Therefore,
\[
    \frac{\#\{(A_1, B_1,A_2, B_2) \in G^4: (A_1, B_1;A_2, B_2) \text{ is a Beauville
    structure}\}}{|G|^4} \geq
\]
\[
    \left(\frac{1}{8} - \eps'_q\right)^2 + \left(\frac{1}{8} - \eps'_q\right)^2
    \geq \frac{1}{32} - \eps_q,
\]
where $\eps_q\rightarrow 0$ as $q \rightarrow \infty$.

Assume that $q$ is odd. In order to obtain an upper bound, observe
that if the orders of $A_1$ and $A_2$ are both even then
Condition~$\textit{(iii)}$ is not satisfied, see~\cite[Lemma
4.2]{Ga}. Denote by $P^e_q$ the probability that a random element in
$G$ has even order, then $P^e_q \geq \frac{1}{4} - \eps'_q$, where
$\eps'_q \rightarrow 0$ as $q \rightarrow \infty$. Indeed, if
$(q-1)/2$ (respectively, $(q+1)/2$) is even, then at least half of
the split (respectively, non-split) elements are of even order,
since any split (respectively, non-split) element belongs to a
cyclic subgroup of order $(q-1)/2$ (respectively, $(q+1)/2$), and at
least half the elements in a cyclic group of even order are of even
order. Therefore,
\[
    \frac{\#\{(A_1, B_1,A_2, B_2) \in G^4: (A_1, B_1;A_2, B_2) \text{ is {\bf not} a Beauville
    structure}\}}{|G|^4} \geq
\]
\[
    \left(\frac{1}{4} - \eps'_q\right)^2 \geq \frac{1}{16} - \eps_q,
\]
where $\eps_q \rightarrow 0$ as $q \rightarrow \infty$.

When $q$ is even the proof is similar, replacing $P^e_q$ by the
probability that a random element in $G$ has order divisible by $3$,
which is at least $\frac{1}{6} - \eps'_q$, where $\eps'_q
\rightarrow 0$ as $q \rightarrow \infty$.

%%%%%%%%%%%%%%%%%%%%%%%%%%%%%%%%%%%%%%%%%%%%%%%%%%%%%%%%%%
\section{Beauville structures for finite simple groups}\label{sect:simple}
%%%%%%%%%%%%%%%%%%%%%%%%%%%%%%%%%%%%%%%%%%%%%%%%%%%%%%%%%%
In this section we briefly describe the probabilistic
group-theoretical methods in proving Theorem~\ref{thm.Beau.FSG} and
Theorem~\ref{thm.Beau.An}. We shall mainly sketch the proof of
Garion, Larsen and Lubotzky~\cite{GLL} that the conjecture of Bauer,
Catanese and Grunewald (Theorem~\ref{thm.Beau.FSG}) holds for
\emph{almost} all finite simple groups of Lie type, as well as
present the proof of Garion and Penegini~\cite{GP} regarding
Beauville structures of alternating groups
(Theorem~\ref{thm.Beau.An}), which is based on the probabilistic
results of Liebeck and Shalev~\cite{LS04}.

As in the probabilistic approach we can ignore finitely many simple
groups, we will not deal here with the sporadic groups, whose
unmixed Beauville structures can be found in~\cite{FMP,GM}. Thus we
shall consider only the alternating groups and the finite simple
groups of Lie type.

Recall that in order to construct an unmixed Beauville structure for
a finite simple (non-abelian) group $G$ one needs to find a
quadruple $(x_1, y_1; x_2, y_2)$ of elements of $G$ satisfying the
three conditions given in Definition~\ref{defn.beau}.

%%%%%%%%%%%%%%%%%%%%%%%%%%%%%%%%%%%%%%%%%%%%%%%%
\subsection{Choosing disjoint conjugacy classes}\label{sect.iii}
%%%%%%%%%%%%%%%%%%%%%%%%%%%%%%%%%%%%%%%%%%%%%%%%
One usually starts by looking for proper conjugacy classes
$X_1,Y_1,Z_1,X_2,Y_2,Z_2$ such that $$\Sigma(x_1,y_1,z_1) \cap
\Sigma(x_2,y_2,z_2)= \{1\}$$ for any $x_i \in X_i,\ y_i \in Y_i,\
z_i \in Z_i$ ($i=1,2$), so that Condition $\textit{(iii)}$ is
satisfied.

For finite simple groups of Lie type, one can choose two maximal
tori $T_1$ and $T_2$, such that if $C_i$ denotes the set of all
conjugates of elements of $T_i$, then $C_1 \cap C_2  = \{1\}$.

For example, let $G=\SL_{r+1}(q)$ ($r>1$), and let $t_1$ and $t_2$
denote elements of $G$ whose characteristic polynomials are
respectively irreducible (of degree $r+1$) and the product of
irreducible polynomials of degree $1$ and $r$. If $T_1$ and $T_2$
denote the centralizers of $t_1$ and $t_2$ respectively, then,
by~\cite[Proposition 7]{GLL}, for all $g \in G$, $T_1 \cap
g^{-1}T_2g = Z(G)$, thus Condition~\textit{(iii)} is satisfied for
$\PSL_{r+1}(q)$.

In fact, for the finite simple groups of Lie type, one can choose
several maximal tori $T_1$ and $T_2$ such that
Condition~\textit{(iii)} is satisfied, see the various choices in
~\cite{FMP,GLL,GM}. However, the number of conjugacy classes of
maximal tori is bounded above by a function of the Lie rank $r$, and
any maximal torus is isomorphic to a product of at most $r$ cyclic
groups, so a similar argument to the one presented in
\S\ref{sect.prob.PSL} implies that the probability that four random
elements generate an unmixed Beauville structure is bounded above by
a function of $r$.

\medskip

One can also choose proper conjugacy classes for the alternating
groups (see~\cite{FMP,FJ,GP,GM}). In~\cite{GP}, Garion and Penegini
used the so-called \emph{almost homogeneous} conjugacy classes
introduced by Liebeck and Shalev~\cite{LS04}.

\begin{defn}\cite{LS04}.
Conjugacy classes in $S_n$ of cycle-shape $(m^k)$, where $n = mk$,
namely, containing $k$ cycles of length $m$ each, are called
\emph{homogeneous}. A conjugacy class having cycle-shape $(m^k
,1^f)$, namely, containing $k$ cycles of length $m$ each and $f$
fixed points, with $f$ bounded, is called \emph{almost homogeneous}.
\end{defn}

By~\cite[Algorithm 3.5]{GP} one can construct for any six integers
$k_1,l_1,m_1$, $k_2,l_2,m_2 \geq 2$, six distinct almost homogeneous
conjugacy classes $X_1,Y_1,Z_1$, $X_2,Y_2,Z_2$ in $A_n$ whose
elements are of orders $k_1,l_1,m_1$, $k_2,l_2,m_2$ respectively,
such that all of them have different numbers of fixed points, thus
satisfying Condition~\textit{(iii)}.

%%%%%%%%%%%%%%%%%%%%%%%%%%%%%%%%%%%%%%%%%%%%%%%%%%%%%%%%%
\subsection{Frobenius formula and Witten's zeta
function}\label{sect.i}
%%%%%%%%%%%%%%%%%%%%%%%%%%%%%%%%%%%%%%%%%%%%%%%%%%%%%%%%%

Condition~\textit{(i)} follows from a classical formula of
Frobenius: If $X$, $Y$ and $Z$ are conjugacy classes in a finite
group $G$, then the number $N_{X,Y,Z}$ of solutions of $xyz = 1$
with $x \in X,$ $y \in Y$ and $z \in Z$ is given by
\begin{equation}\label{eq.frob}
N_{X,Y,Z} = \frac{|X|\cdot |Y| \cdot |Z|}{|G|} \sum_{\chi \in
\Irr(G)} \frac{\chi(x)\chi(y)\chi(z)} {\chi(1)},
\end{equation}
where $\Irr(G)$ denotes the set of complex irreducible characters of
$G$.

Usually, the main contribution to this character sum comes from the
trivial character, and the absolute sum on all other characters is
negligible. In order to show this, one needs to bound the
absolute value $|\chi(g)|$ for an irreducible character $\chi$ and
an element $g$, from any of the above conjugacy classes
(see~\S\ref{sect.char} for details). If this value can be
effectively bounded then one deduces that the value of the sum of
the contribution of all \emph{non-trivial} characters to
\eqref{eq.frob} is bounded above by some global constant (depending
only on the sizes of $G$ and the conjugacy classes) multiplied by
the sum $\sum_{\chi \neq 1} \chi(1)^{-1}$. So it remains to prove
that the letter sum converges to $0$ as $|G| \rightarrow \infty$.

Therefore, a key role in the probabilistic approach is played by the
so called \emph{Witten zeta function}, which is defined by
$$\zeta^G(s) = \sum_{\chi \in \Irr(G)} \chi(1)^{-s}.$$ It was
originally defined and studied by Witten~\cite{Wit} for Lie groups.
For finite simple groups it was studied and applied in detail by
Liebeck and Shalev~\cite{LS04,LS05C,LS05F}, who proved the following
desired results.

\begin{thm}\cite[Theorem 1.1]{LS05F}, \cite[Corollary 1.3]{LS05C} and \cite[Corollary 2.7]{LS04}.
\label{thm.wit} \\
Let $G$ be a finite simple group.

\begin{itemize}
\item If $s > 1$ then $\zeta^G(s) \rightarrow 1$ as $|G| \rightarrow
\infty$.
\item If $s > 2/3$ and $G \neq \PSL_2(q)$ then $\zeta^G(s) \rightarrow
1$ as $|G| \rightarrow \infty$.
\item If $s > 1/2$ and $G \neq \PSL_2(q),\PSL_3(q),\PSU_3(q)$
then $\zeta^G(s) \rightarrow 1$ as $|G| \rightarrow \infty$.
\item If $s > 0$ and $G = A_n$ then $\zeta^G(s) \rightarrow 1$ as $|G|
\rightarrow \infty$. Moreover, $\zeta^G(s) = 1 + O(n^{-s}).$
\end{itemize}
\end{thm}

An alternative approach to prove Condition~\textit{(i)} for the
finite simple groups of Lie type, which was successfully applied
in~\cite{FMP,GM} and~\cite[\S4]{GLL}, is based on the following
result of Gow~\cite{Gow}: if $X$ and $Y$ are conjugacy classes of
regular semisimple elements in a finite Lie type group $G$, then the
set $XY$ contains every non-central semisimple element of $G$.

%%%%%%%%%%%%%%%%%%%%%%%%%%%%%%%%%%%%%%%%%%%%%%%%%%%%%%%%
\subsection{Character estimates in finite simple groups}\label{sect.char}
%%%%%%%%%%%%%%%%%%%%%%%%%%%%%%%%%%%%%%%%%%%%%%%%%%%%%%%%

A crucial part in the proof is to estimate character values in
finite simple groups. More precisely, one needs to bound the
absolute value $|\chi(g)|$ for an irreducible character $\chi$ and
an element $g$, from any of the conjugacy classes chosen
in~\S\ref{sect.iii}. We therefore recall some useful results for the
finite simple groups of Lie type and the alternating groups.

\begin{lemma}\cite[Corollary 4]{GLL}.
Let $\chi$ be an irreducible character of $G = \SL_{r+1}(q)$.

\begin{itemize}
\item If $t_1\in G$ has an irreducible characteristic polynomial, then
$| \chi(t_1)|  \le \frac{2(r+1)^2}{\sqrt3}.$

\item If $r\ge 2$ and the characteristic polynomial of $t_2\in G$ has
an irreducible factor of degree $r$,  then $|\chi(t_2)|  \le
\frac{2r^2}{\sqrt3}.$
\end{itemize}
\end{lemma}

More generally, by~\cite[Proposition 7]{GLL}, there exist an
absolute constant $c$ such that for every sufficiently large group
of Lie type $G$ (of Lie rank $r$), there exist maximal tori $T_1$
and $T_2$ as in~\S\ref{sect.iii}, and for every regular $t \in
T_1\cup T_2$ and every irreducible character $\chi$ of $G$,
$$|\chi(t)| \leq cr^3.$$

\begin{lemma}\cite[Proposition 2.12]{LS04}.
Let $\pi \in S_n$ have cycle-shape $(m^k,1^f)$. Then for any $\chi
\in \Irr(S_n)$ we have
$$ |\chi(\pi)| \leq c \cdot (2n)^{\frac{f+1}{2}}\chi(1)^{\frac{1}{m}},$$
where $c$ depends only on $m$.
\end{lemma}

It is also interesting to provide an upper bound on the \emph{character ratio}
$|\chi(g)/\chi(1)|$, where $G$ is a finite simple group, $g \in G$ and $\chi$ is an
irreducible character. Gluck and Magaard~\cite{GlM} computed these bounds for the 
finite classical groups.
Such bounds play a crucial role in the proof of a conjecture of Guralnick and Thompson~\cite{GT}, which is related to the inverse Galois problem, namely, 
which finite groups occur as Galois groups of algebraic number fields 
(for details see~\cite{FM}). 

%%%%%%%%%%%%%%%%%%%%%%%%%%%%%%%%%%%%%
\subsection{Finding generating pairs}\label{sect.ii}
%%%%%%%%%%%%%%%%%%%%%%%%%%%%%%%%%%%%%
In order to prove Condition~\textit{(ii)} one should show that the
set of solutions of $xyz = 1$ with $x \in X,$ $y \in Y$ and $z \in
Z$ contains a generating pair of $G$, namely, one should avoid pairs
$(x,y)$ contained in maximal subgroups of $G$. Probabilistically,
one expects such a result to hold since almost all pairs of elements
in a finite simple (non-abelian) group are generating pairs (see
\cite{Dix,KL,LS95}).

Namely, one should estimate the sum
\[
\label{eq.non.max} \sum_{M\in \max G}|\{(x,y,z)\mid \\ x\in X\cap M,
y \in Y\cap M, z \in Z \cap M,  xyz=1 \}|,
\]
where $\max G$ denotes the set of maximal proper subgroups of $G$.
This quantity is bounded above by
\[\sum_{M\in \max G} |M|^2 = |G|^2
\sum_{M\in \max G} \frac 1{[G:M]^2} \le
\frac{|G|^2}{m(G)^{1/2}}\sum_{M\in \max G}[G:M]^{-3/2},
\]
where $m(G)$ is the minimal index of a proper subgroup of $G$ or,
equivalently, the minimal degree of a non-trivial permutation
representation of $G$. By estimates of Landazuri and
Seitz~\cite{LS}, there exists an absolute constant $c$ such that if
$G$ is a finite simple group of Lie type $G$ of Lie rank $r$ then
$m(G) \ge c q^r$. Hence, again one should estimate a 'zeta function'
encoding the indices of maximal subgroups of finite simple groups of
Lie type, which was investigated by Liebeck, Martin and
Shalev~\cite{LMS}, who proved the following desired result.

\begin{thm}\cite{LMS}.
If $G$ is a finite simple group, and $s > 1$, then
$$\lim_{|G|\to\infty}\sum_{M\in \max G}[G:M]^{-s}\to 0.$$
\end{thm}

Alternatively, in~\cite{FMP,GM} they relied on more delicate results
about maximal subgroups in finite simple groups of Lie type
containing special elements, called \emph{primitive prime divisors},
of Guralnick, Pentilla, Praeger and Saxl~\cite{GPPS}.

\medskip

For sufficiently large alternating groups, Conditions~\textit{(i)}
and~\textit{(ii)} follow from the following result of Liebeck and
Shalev~\cite{LS04}. If $(k,l,m)$ is hyperbolic then the probability
that three random elements $x,y,z\in A_n$, with product $1$, from
almost homogeneous classes $X,Y,Z$, of orders $k,l,m$ will generate
$A_n$, tends to 1 as $n \rightarrow \infty$. Moreover, this
probability is $1-O(n^{-\mu})$, where $\mu=1-(1/k+1/l+1/m)$.

%%%%%%%%%%%%%%%%%%%%%%%%%%%%%
\subsection*{Acknowledgement}
%%%%%%%%%%%%%%%%%%%%%%%%%%%%%
The author would like to thank her co-organizers, Ingrid Bauer and
Alina Vdovina, and all the participants in the workshop ``Beauville
surfaces and groups 2012'' for their assistance and useful conversations, 
as well as the University of Newcastle for hosting the workshop.
The author is grateful to the referee for pointing out further relevant references.

%%%%%%%%%%%%%%%%%%%%%%%%%%%%%%%%%%%%%%%%%%%%%%%%%%%%%%%
%%%%%%%%%%%%%%%%%%%%%%%%%%%%%%%%%%%%%%%%%%%%%%%%%%%%%%%

\end{document}